\documentstyle{article}
\input epsf.tex
\begin{document}

\title
{
About a conjecture for uniformly isochronous polynomial centers
}

\author{
Evgenii P. Volokitin%
\thanks{Partially supported by Grant 05--01--00302 from the Russian
Foundation of Basic Research.}\\
Sobolev Institute of Mathematics, Novosibirsk, 630090, Russia\\
e-mail: volok@math.nsc.ru
}
\date{}
\maketitle

\begin{abstract}
We study a specific family of uniformly isochronous polynomial systems.
Our results permit to solve a problem about centers of such systems.
\end{abstract}

Classification: Primary 34C05; Secondary 34C25


\quad

{\bf 1.}
Consider the planar autonomous system of ordinary differential
equations
\refstepcounter{equation}
\label{1}
$$
\begin{array} {ll}
\dot x = -y + x H (x,y),\\
\dot y =  x +y H (x,y),\\
\end{array}
\eqno{(\ref{1})}
$$
where $H(x,y)$ is a polynomial in $x$ and $y$ of degree $n$,
and $H(0,0)=0$. This system has only one singular point at $O(0,0)$
which is the center of the linear part of the system. The solutions
of this system move around the origin with constant angular speed,
and the origin is so a uniformly isochronous singular point.

The problem of characterizing uniformly isochronous centers has attracted
attention of several authors; see \cite{1}--\cite{4} and references therein.
In particular, the following problem was posed:

{\it
It is true that all centers for uniformly isochronous polynomial systems
are either reversible or admit a nontrivial polynomial commuting system?
}

The problem first appeared in \cite{2} and it was mentioned as one of the open
questions in \cite{3}. It was also marked by the reviewers in
Zbl. Math. 1037.34024 and in MR1963468 (2004b:34090).
We prove the following proposition which permits to give a negative answer
to the problem.

{\bf Theorem 1.} Let a uniformly isochronous polynomial system has the form
\refstepcounter{equation}
\label{2}
$$
\begin{array} {ll}
\dot x = -y + x Q (x,y) \sum_{i=0}^{m} a_i (x^2+y^2)^i,\\
\dot y =  x +y Q (x,y) \sum_{i=0}^{m} a_i (x^2+y^2)^i,
\end{array}
\eqno{(\ref{2})}
$$
where $Q(x,y)$ is a homogeneous polynomial in $x,y$ of degree $k$ and
\refstepcounter{equation}
\label{3}
$$
\int_0^{2\pi} Q(\cos \vartheta, \sin \vartheta) d\vartheta=0.
\eqno{(\ref{3})}
$$
Then the origin is a center of (\ref{2}). The center is of type $B^\nu$
with $\nu \leq k$, and a <<generic>> center is of type $B^1$ when $k$ is
odd or of type $B^2$ when $k$ is even.

{\bf Proof.}
System (\ref{2}) can be written as a single separable equation
\refstepcounter{equation}
\label{4}
$$
\frac{d \varrho}{d \vartheta}=
\varrho^{k+1} Q(\cos \vartheta, \sin \vartheta) R(\varrho)
\eqno{(\ref{4})}
$$
with $\varrho, \vartheta$ polar coordinates and $R(\varrho)=
\sum_{i=0}^{m} a_i \varrho^{2i}$.

Equation (\ref{4}) has a solution $\varrho \equiv 0$ which is defined
for all $\vartheta$. Therefore every solution $\varrho(\vartheta)$ with the
initial value $\varrho(0)=\varrho_0$ where $\varrho_0>0$ is small enough
is defined for $\vartheta \in [0, 2\pi]$ and satisfies the condition
\refstepcounter{equation}
\label{5}
$$
\int_0^{\vartheta} Q(\cos \varphi, \sin \varphi) d\varphi=
\int_{\varrho_0}^{\varrho(\vartheta)} \frac{dr}{r^{k+1} R(r)}.
\eqno{(\ref{5})}
$$

From (\ref{3}) we conclude that the solution is $2\pi$-periodic,
so that the origin is a center.
The first part of the theorem is proved.

By \cite{5}, the center of (\ref{2}) is of type $B^{\nu}$, and the boundary
of the center region is the union of $\nu$ open unbounded trajectories
with $\nu \leq n=k+2m$.

The circles $x^2+y^2=\varrho_i^2$ where $\varrho_i$ are roots
of the equations $R(\varrho)=0$ are trajectories of (\ref{2}).
All of them lie in the center region.
$$
\epsfbox{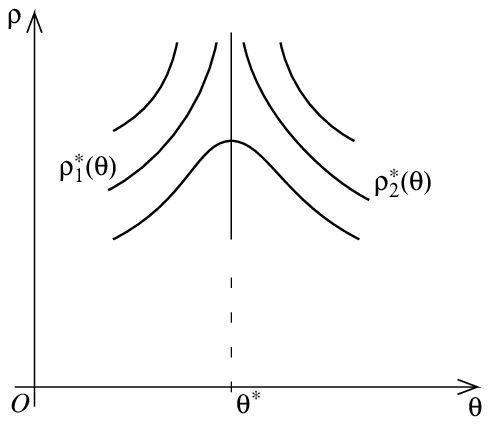} \qquad \qquad \qquad
\epsfbox{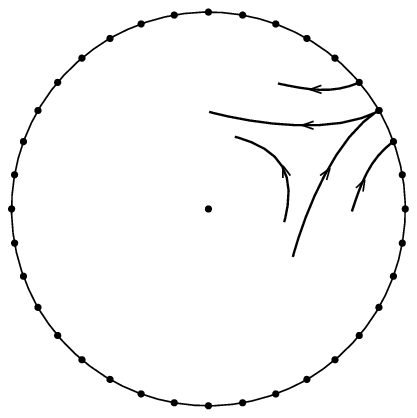}
$$

\hspace{2cm} Fig. 1 \hspace{5.5cm} Fig. 2

Unbounded trajectories of (\ref{2}) correspond to unbounded solutions
of (\ref{4}). Studing the behaviour of solution curves of (\ref{4}) for
large $\varrho$, it can be shown that for every null isocline
$\vartheta=\vartheta^*$ where solutions have a maximum there exist two
solutions $\varrho_1^*(\vartheta), \varrho_2^*(\vartheta)$,
for which the isocline is a vertical asymptote
$$
\lim_{\vartheta \rightarrow \vartheta^*-0} \varrho_1^*(\vartheta)=+\infty,
\lim_{\vartheta \rightarrow \vartheta^*+0} \varrho_2^*(\vartheta)=+\infty,
$$
(see Fig. 1).

In this situation, there is a relevant equilibrium point at infinity
in the intersection of the equator of the Poincar\'e sphere with
the ray $x=\varrho \cos \vartheta^*, y=\varrho \sin \vartheta^*, \varrho>0$.
The point has one hyperbolic sector with two separtrices
which correspond to the solutions
$\varrho_1^*(\vartheta), \varrho_2^*(\vartheta)$ (see Fig. 2).
The boundary of the center region consists of such separatrices.
The number $\nu$ of these equilibrium points coinsides with the
number of the null isoclines of direction field (\ref{4}) where
solutions has a maximum for large $\varrho$. These isoclines are
vertical lines $\vartheta=\vartheta^*_i$, where the values of
$\vartheta^*_i$ are determined from the conditions
$Q(\cos\vartheta,\sin\vartheta)=0, 0\leq \vartheta <2\pi$.
Hence we have the estimate $\nu\leq k$ when we describe type $B^{\nu}$
of the center of system (\ref{2})%
\footnote{%
If $k$ is even our trigonometric polynomial
$Q(\cos\vartheta,\sin\vartheta)$ has a period equal to $\pi$ (but not
$2\pi$ as it takes place for odd $k$). Therefore (\ref{4})
has an even number of the blocks discussed above and the relevant
equilibrium points are disposed at the diameters of the Poincar\'e
sphere. It may be noted that system (\ref{2}) is $O$-symmetric in this
case.}.
The upper bound $k$ can be attained by $\nu$. As an example we can
consider system (\ref{2}) with
$Q(\cos\vartheta,\sin\vartheta)=\sin k \vartheta$ and $a_i$ arbitrary
real numbers.

In a <<generic>> situation, (\ref{4}) has no solution for which two different
null isoclines are its asimptotes. Therefore in such a situation
the solution curve which separates bounded and unbounded solutions
has a minimum number of discontinuity points within $[0, 2\pi]$,
it has one point when $k$ is odd, and it has two points when $k$ is even.

Hence a <<generic>> center is of type $B^1$ when $k$ is odd
or type $B^2$ when $k$ is even (see Fig. 3).
$$
\epsfbox{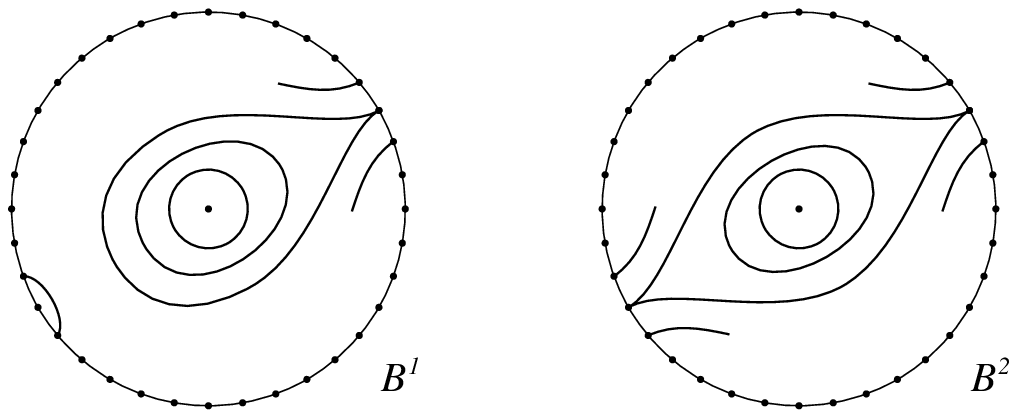}
$$
\begin{center}
{Fig. 3}
\end{center}

The theorem is proved.

{\bf Remark.} Theorem 2.1 from \cite{5} about centers of homogeneous systems
is a particular case of our theorem for $m=0, a_0=1$.

The functions
$$
f_1(x,y)=x^2+y^2, f_2(x,y)=\sum_{i=0}^{m} a_i (x^2+y^2)^i
$$
are invariants for (\ref{2}) with respective cofactors
$$
\begin{array}{ll}
K_1(x,y)=2Q(x,y)\sum_{i=0}^{m} a_i (x^2+y^2)^i,\\
K_2(x,y)=2Q(x,y)\sum_{i=0}^{m} i a_i (x^2+y^2)^i.
\end{array}
$$

We have
$$
\frac{k+2}{2} K_1(x,y) + K_2(x,y)= div,
$$
where
$div$ is the divergence of (\ref{2}). In this case the function
$\mu(x,y)=f_1^{(k+2)/2}f_2$ is a reciprocal integrating factor
of Darboux form\footnote%
{About algebraic invariants and Darboux's method of integration see
\cite{7}, \cite{8}, for example.}.
The factor gives information about our system. For instance, it may be used
to find a first Darboux integral of (\ref{2}) \cite{6}.
A first integral of (\ref{2}) may be found from (\ref{5}) also.

It is obvious that (\ref{2}) commutes with the system
\refstepcounter{equation}
\label{6}
$$
\begin{array}{ll}
\dot x = x (x^2+y^2)^{k/2} \sum_{i=0}^{m} a_i (x^2+y^2)^i,\\
\dot y = y (x^2+y^2)^{k/2} \sum_{i=0}^{m} a_i (x^2+y^2)^i.
\end{array}
\eqno{(\ref{6})}
$$

If $k$ is even (\ref{6}) gives a polynomial commuting system without a
linear part. If $k$ is odd we have a non-polynomial commuting system.
Nevertheless a polinomial commuting system may exist in the case of
odd $k$. For example, if (\ref{2}) is homogeneous ($m=0, a_0=1$)
then there exists a polynomial system which commmutes with (\ref{2}) \cite{9}.

Using Theorem 1, we may consruct an example of an uniformly isochronous
system which is not reversible and commutes with no polynomial system.

If a system is reversible then its trajectories are symmetric with respect
a common symmetric line%
\footnote{%
Necessary and sufficient conditions for reversibility of planar
analytic vector fields were derived in \cite{10}.}.
If a symmetric line of system (\ref{2}) is
$x \sin \vartheta^{*} - y \cos \vartheta^{*}=0$,
then the vertical line $\vartheta=\vartheta^{*}$ is the symmetric axis of
the graph of the trygonometric polynomial $Q(\cos \vartheta, \sin \vartheta)$
and it is the symmetric axis of solution curves of vector field (\ref{4}).

The problem about conditions for the existence of polynomial commuting
systems for uniformly isochronous polynomial systems was considered in
\cite{2}, \cite{3}. In partucular, it was proved that system (\ref{1})
commutes with a polynomial system if anf only if the function $H(x,y)$
satisfies  one of the following two conditions:

\noindent 1)
\refstepcounter{equation}
\label{7}
$$
H(x,y)=P_{2l}(x,y)\sum_{j=0}^r a_j (x^2+y^2)^j
\eqno(\ref{7})
$$
where $P_{2l}(x,y)$ is a homogeneous polynomial of degree $2l, l\geq 0$;

\noindent 2) there are homogeneos polynomials $\alpha_l, \beta_l$ of
order $l$
($l\leq n$, $l$ divides $n$), verifying
$x \partial_y \beta_l-y \partial_x \beta_l=l \alpha_l$, such that
\refstepcounter{equation}
\label{8}
$$
H(x,y)=\alpha_l \sum_{k=0}^{n/l-1} a_k \beta_l^k.
\eqno(\ref{8})
$$

So, to consruct our example it is sufficient to take a system
of the form (\ref{2}) where the homogeneous polynomial $Q(x,y)$ is of
an odd degree\footnote%
{In this case (\ref{3}) is fulfilled and the function $H(x,y)$
is not of the form (\ref{7}).},
the graph of the trygonometric polynomial $Q(\cos \vartheta, \sin\vartheta)$
has no symmetric axes and the numbers $m, a_i$ are such that the function
$$
H(x,y)=Q(x,y) \sum_{i=0}^{m} a_i (x^2+y^2)^i
$$
is not of the form (\ref{8}).

Let us set
$$
Q(x,y)=y^3-3x y^2+2x^2 y=y(x-y)(2x-y), m=1, a_0=a_1=1.
$$
Then we have
\refstepcounter{equation}
\label{9}
$$
\begin{array}{ll}
\dot x=-y+x (y^3-3x y^2+2 x^2 y)(1+x^2+y^2),\\
\dot y=x+y (y^3-3x y^2+2 x^2 y)(1+x^2+y^2).
\end{array}
\eqno(\ref{9})
$$

According to Theorem 1, (\ref{9}) has a center (isochronous) at the origin.

The function
$$
I(x,y)=\frac{r^6}{(1-3r^2-4x^3-3x y^2-3y^3-3r^3\arctan r)^2}, r^2=x^2+y^2.
$$
is a first integral of (\ref{9}) obtained from (\ref{5}).

It is evident that the graph of $Q(\cos \vartheta, \sin\vartheta )$
has no symmetric axes%
\footnote{%
It may be shown that if the graph of the homogeneous trygonometric
polynomial of degree~3
$$
T_3(\vartheta)=a_1 \cos \vartheta+a_3 \cos 3\vartheta+b_1 \sin \vartheta+
b_3 \sin 3\vartheta
$$
has a simmetric axes then its coefficients satisfy the condition
$$
a_1 b_3(a_1^2-3b_1^2)=a_3 b_1(3a_1-b_1^2).
$$
},
and therefore system (\ref{9}) is non-reversible.

It easy to veryfy that system (\ref{9}) may fail to commute with any
nonprortional polynomial systems. This fact follows from the impossibility
to present the function
$$
H(x,y)=(y^3-3x y^2+2 x^2 y)(1+x^2+y^2)\equiv H_3(x,y)+H_5(x,y)
$$
in the form (\ref{8}) but we can also prove it immediately.

Indeed, let system (\ref{9}) commutes with a polynomial system of degree $n$
\refstepcounter{equation}
\label{10}
$$
\begin{array}{ll}
\dot x= R(x,y)\equiv R_1(x,y)+R_2(x,y)+\dots+R_n(x,y),\\
\dot y= S(x,y)\equiv S_1(x,y)+S_2(x,y)+\dots+S_n(x,y),
\end{array}
\eqno{(\ref{10})}
$$
where $R_i(x,y), S_i(x,y)$ are homogeneous polynomials of degree $i$.

Then the Lie bracket between vector fields (\ref{9}), (\ref{10})
is equal to zero:
$$
[(-y+x H(x,y), x+y H(x,y))^T, (R(x,y), S(x,y))^T]=(0,0)^T.
$$

In particular, we have that terms of highest degree are equal to zero:
$$
[(x H_5(x,y), y H_5(x,y))^T, (R_n(x,y), S_n(x,y))^T]=(0,0)^T.
$$

After transformations taking into account Euler's theorem for homogeneous
functions this equality may be written in the form
$$
\begin{array}{ll}
(x H_{5x}(x,y)+(1-n)H_5(x,y))R_n(x,y)+x H_{5y}(x,y) S_n(x,y)=0,\\
y H_{5x}(x,y) R_n(x,y)+(y H_{5y}(x,y)+(1-n)H_{5}(x,y)) S_n(x,y)=0.
\end{array}
$$

The linear system for determing the polynomials $R_n(x,y), S_n(x,y)$
has a nontrivial solution if its determinant $\Delta$ is equal zero:
$$
\begin{array}{ll}
\Delta\equiv& (x H_{5x}(x,y)+(1-n)H_5(x,y)) (y H_{5y}(x,y)+(1-n)H_{5}(x,y))-\\
 &x y H_{5x}(x,y) H_{5y}(x,y)=0,
\end{array}
$$

Taking into account the fact that
$x H_{5x}(x,y) + y H_{5y}(x,y) =5 H_{5}(x,y)$, we have
$$
\Delta=(1-n)(6-n)H_{5}^2 (x,y)=0.
$$

Therefore we must have that $n=6$ or $n=1$. Straightforward calculations
using the software package {\it Mathematica} show that in this case the
commuting system (\ref{10}) is proportional to system (\ref{9})%
\footnote{According with (\ref{6}) system (\ref{9}) commutes with the system
$$
\dot x=x(x^2+y^2)\sqrt{x^2+y^2}(1+x^2+y^2),
\dot y=y(x^2+y^2)\sqrt{x^2+y^2}(1+x^2+y^2).\\
$$
}.

Hence system (\ref{9}) has a center but it is non-reversible and it
commutes with no polynomial system nonproportional to it.

We derive that the answer to the question from \cite{2}, \cite{3} is negative.

{\bf 2.} We can generalize the first part of Theorem 1.

{\bf Theorem 2.} Let the polynomial $H(x,y)$ in (\ref{1}) has the form
$$
H(x,y)=q(x,y) h(x^2+y^2, p(x,y)),
$$
where $h(u,v)$ is a polynimial, $p(x,y), q(x,y)$ are homogeneous polynomials
of the same degree $k$ and $q(x,y)=c(x p_{y}(x,y)-y p_{x}(x,y))$.
Then the origin is a center of (\ref{1}).

{\bf Proof.}

In the case under study system (\ref{1}) can be written as a single
equation of the form
$$
\frac{d \varrho}{d \vartheta}=
c \varrho^{k+1}h(\varrho^2,\varrho^{k} f(\vartheta)) f'(\vartheta)
$$
with $\varrho, \vartheta$ polar coordinates and
$f(\vartheta)=p(\cos\vartheta, \sin\vartheta)$.

It is clear that solutions of this equation are functions of $f(\vartheta)$.
The function $f(\vartheta)$ is 2$\pi$-periodic. Then solutions with initial
values which are small enough are 2$\pi$-periodic functions also.
So, the origin is a center.


\begin{thebibliography}{99}

\bibitem{1}
J. Chavarriga, M. Sabatini;
{\it A survey of isochronous centers},
Qualitative Theory of Dynamical Systems. 1999. Vol. 1. No. 1. P. 1--70.

\bibitem{2}
A. Algaba, M. Reyes, and A. Bravo;
{\it Geometry of the uniformly isochronous centers with polynomial
commutators},
Differential Equations Dynam. Systems. 2002. Vol. 10. No. 3-4. P. 257--275.

\bibitem{3}
A. Algaba, M. Reyes;
{\it Centers with degenerate infinity and their commutators},
J. Math. Anal. Appl. 2003. Vol. 78. No. 1. P. 109--124.

\bibitem{4}
A. Algaba, M. Reyes;
{\it Computing center conditions for vector fields with constant angular
speed},
J. Comput. Appl. Math. 2003. Vol. 154. No. 1. P. 143--159.

\bibitem{5}
R. Conti,
{\it Uniformly isochronous centers of polynomial systems in $R^2$},
Elworthy K. D. (ed.) et al., Differential equations, dynamical
systems, and control science. New York: Marcel Dekker.
Lect. Notes Pure and Appl. Math. 152, 21--31 (1994).

\bibitem{6}
J. Chavarriga, H. Giacomini, and J. Gin\'e;
{\it The null divirgence factor},
Publ. Mat., Barc. 1997. Vol. 41. No. 1. P. 41--56.

\bibitem{7}
J. M. Pearson, N. G. Lloyd, and C.J. Christopher;
{\it Algorithmic derivation of centre conditions},
SIAM Review. 1996. Vol. 38. No. 4. P. 619--636.

\bibitem{8}
D. Schlomiuk,
{\it Algebraic and geometric aspects of the theory of polynomial vector
fields},
Schlomiuk, Dana (ed.), Bifurcations and periodic orbits of vector fields.
Proceedings of the NATO Advanced Study Institute and S\'eminaire de
Math\'ematiques Sup\'erieures, Montr\'eal, Canada, July 13-24, 1992.
Dordrecht: Kluwer Academic Publishers. NATO ASI Ser., Ser. C, Math. Phys. Sci.
408, 429-467 (1993).

\bibitem{9}
L. Mazzi, M. Sabatini;
{\it Commutators and linearizations of isochronous centers},
Atti Acad. Naz. Lincei Cl. Sci. Fis. Mat. Natur. Rend. Lincei (9)
Mat. Appl. 2000. Vol. 11. No 2. P. 81--98.

\bibitem{10}
C. B. Collins,
{\it Poincar\'e's reversibility conditions},
J. Math. Anal. Appl. 2001. Vol. 259. No. 1. P. 168--187.

\end{thebibliography}
\end{document}